\documentclass[amsfonts]{article}
\usepackage{amsmath,dsfont}
\usepackage{amssymb}
\newtheorem{teo}{Theorem}[section]
\newtheorem{prop}[teo]{Proposition}
\newtheorem{lm}[teo]{Lemma}

\newtheorem{rem}[teo]{Remark}

\newcommand{\RR}{{\mathbb{R}}}
\newcommand{\CC}{{\mathbb C}}
\newcommand{\NN}{{\mathbb N}}

\newcommand{\der}{\partial}
\newcommand{\ep}{\varepsilon}
\newcommand{\ov}{\overline}

\newcommand{\om} {\Omega}

\newcommand{\Gm}{\mathbf{G}}

\def\qed{\hfill$\square$\vspace{0.5cm}}    

\begin{document}

\title{Lipschitz stability for the electrical impedance tomography problem: the complex case.}

\author{Elena Beretta\thanks{ Dipartimento di Matematica ``G. Castelnuovo''
Universit\`a
 di Roma ``La Sapienza'', Piazzale Aldo Moro 5, 00185 Roma, Italy
(Email: beretta@mat.uniroma1.it).} \and Elisa Francini\thanks{Dipartimento di Matematica "U. Dini",
Viale Morgagni 67A, 50134 Firenze, Italy
(Email: francini@math.unifi.it).} }
\date{}
\maketitle
\begin{abstract}
In this paper we investigate the boundary value problem
\[\label{ab}
    \left\{\begin{array}{rcl}
             \mbox{div}\left(\gamma\nabla u\right)&= & 0\quad\mbox{in}\quad\om \\
             u&= &f  \quad\mbox{on}\quad\der\om,\\
           \end{array}
    \right.
\]
where $\gamma$ is a complex valued $L^\infty$ coefficient, satisfying a strong ellipticity condition.
In Electrical Impedance Tomography, $\gamma$ represents the admittance of a conducting body. An interesting issue is the one of determining $\gamma$  uniquely and in a stable way from the knowledge of the Dirichlet-to-Neumann map $\Lambda_\gamma$. Under the above general assumptions this problem is an open issue.

In this paper we prove that, if we assume a priori that $\gamma$ is piecewise constant with a bounded known number  of unknown values, then Lipschitz continuity of $\gamma$ from $\Lambda_\gamma$ holds.
\end{abstract}

\section{Introduction}\label{sec1}
In this paper we investigate stability for the inverse problem of electrical impe\-dance tomography.
More precisely we consider the following problem: let $u\in H^1(\om)$ be the solution to
\begin{equation}\label{1.1}
    \left\{\begin{array}{rcl}
             \mbox{div}\left(\gamma\nabla u\right)&= & 0\quad\mbox{in}\quad\om \\
             u&= &f  \quad\mbox{on}\quad\der\om,\\
           \end{array}
    \right.
\end{equation}
where $\om\subset\RR^n$, $n\geq 2$ is a bounded connected domain, $\gamma$ is a complex valued function representing the admittivity coefficient, it is bounded and satisfies the ellipticity condition $\Re\gamma\geq\lambda^{-1}>0$ a.e. in $\om$ and $f\in H^{1/2}(\der\om)$.

\noindent The Dirichlet to Neumann map $\Lambda_\gamma$ is the operator $\Lambda_\gamma:H^{1/2}(\der\om)\rightarrow H^{-1/2}(\der\om)$ given by
\[
\Lambda_\gamma f=\gamma{\frac{\der u}{\der\nu}}_{|_{\der\om}},
\]
where $\nu$ is the exterior unit normal vector to $\der\om$.

The mathematical formulation of the inverse problem of impedance tomography is to determine the admittivity $\gamma$ from the knowledge
of the Dirichlet to Neumann map $\Lambda_\gamma$.

This problem has several important applications in fields like medical imaging and nondestructive testing of materials.
We refer to the review papers by Borcea (\cite{Bo}) and to (\cite{CIN}) for a wide bibliography on relevant examples of applications.
We want to point out that equation (\ref{1.1})  also appears  in the study of a model for electrical conduction in biological tissues
as the asymptotic limit of an elliptic equation with memory when subjected to periodic Dirichlet boundary conditions (see \cite{AABG}
and  \cite {AABG2}).

 For $n\geq 3$ the uniqueness result by Sylvester and Uhlmann (\cite{SU}) obtained for real conductivities applies also to the
complex case (cfr. \cite{Bo}). For $n=2$ the first contribution on the unique determination of $\gamma$ by $\Lambda_\gamma$ was
given in \cite{F} where the author proved uniqueness provided the imaginary part of $\gamma$ is sufficiently small.
In 2008 Bukhgeim in \cite{Bu} generalized this result to an arbitrary sufficiently smooth admittance.

The problem of determining uniquely an arbitrary $L^\infty$ admittivity from the Dirichlet to Neumann map is completely open
even in the real case when $n\geq 3$. The only result known is the one of Astala and Paivarinta who proved uniqueness of
real $L^\infty$ conductivities from the Dirichlet to Neumann map in the two dimensional case, (cf. [AP]).

The main topic of our paper is to investigate continuous dependence of $\gamma$ on $\Lambda_\gamma$ when the admittivity is
an $L^\infty$ function of a particular form.

In general for arbitrary conductivities  it  is well known that this problem is severely ill-posed. If $\gamma$ is a real
valued coefficient satisfying suitable a-priori smoothness assumptions, Alessandrini proved in \cite{A} a log-type stability
estimate for $n\geq 3$; the same type of stability was proved in \cite{BFR} for $n=2$ for  H\"{o}lder continuous conductivities.
 Such estimates are optimal (see \cite{M}).
Clearly one expects the same kind of ill-posedness also in the complex case.

On the other hand, in many applications one has to disposal additional a priori information on the unknown function the might lead to better stability bounds. In \cite{AV} Alessandrini and Vessella assume that the real conductivity $\gamma$ is of the form
\begin{equation}
\gamma(x)=\sum_{j=1}^N \gamma_j \mathds{1}_{D_j}(x)\quad \mbox{a.e. in }\om,
\end{equation}
where $D_j$ are known disjoint Lipschitz domain and $\gamma_j$ are unknown real numbers.

Assuming ellipticity and $C^{1,\alpha}$ regularity at the interfaces joining contiguous domains $D_j$  and at $\der\om$ they prove
Lipschitz continuous dependence of $\gamma$ on $\Lambda_\gamma$.
The key ingredients in their proof are, on one hand the use of the Green's function and its asymptotic behaviour near the regular interfaces, on the other hand the use of global $C^{\alpha}$ regularity estimates of solutions and local $C^{1,\alpha}$ regularity estimates in a neighborhood of the smooth interfaces.

In this paper we generalize the result in \cite{AV} to the complex equation (\ref{1.1}).  More precisely we show that if $\gamma^{(1)}$ and $\gamma^{(2)}$ are of the form
\[
\gamma^{(k)}(x)=\sum_{j=1}^N \gamma_j^{(k)}\mathds{1}_{D_j}(x),\quad k=1,2
\]
with $\Re \gamma^{(k)}\geq \lambda^{-1}>0$ for $k=1,2$ and assuming  that the interfaces joining contiguous domains contain a flat  portion then
\[
\|\gamma^{(1)}-\gamma^{(2)}\|_{L^\infty(\om)}\leq C \|\Lambda_{\gamma^{(1)}}-\Lambda_{\gamma^{(2)}}\|_{\mathcal{L}\left(H^{1/2}(\partial\Omega),H^{-1/2}(\partial\Omega)\right)},\]
where $C$ depends on $\om$, $\lambda$ and $N$ and diverge to $+\infty$ exponentially as $N\to+\infty$.

Our approach follows the one of Alessandrini and Vessella of constructing singular solutions and  of studying their asymptotic behaviour when the singularity approaches the discontinuity interface.

Observe that if $\Re\gamma\geq\lambda^{-1}>0$  the complex equation (\ref{1.1}) is equivalent to a two by two strongly elliptic system with $L^{\infty}$ coefficients.
One relevant difference with the conductivity case treated in \cite{AV} is that  in the case of real $L^{\infty}$ conductivities existence of the Green's function in $\Omega$ is guaranteed by the results contained in \cite{LSW}  while
for equation (\ref{1.1}), $L^{\infty}$ admittivities and $n\geq 3$ the existence of the Green's function in the whole domain $\Omega$ is not known due to the lack of a maximum principle and of De Giorgi-Nash type regularity estimates for this type of equations.

We are able to bypass these difficulties observing that in order to derive our result  it is enough to construct and to study the behaviour of singular solutions in a Lipschitz subset $K$ (defined in Section 3) of a slightly enlarged domain $\Omega_0$  containing the smooth portion of the interfaces and determining its  asymptotic behaviour near  the interfaces. On the other hand in the domain $K$, using the estimates for elliptic systems obtained by Li and Nirenberg in  \cite{LN},  we have that solutions to equation (\ref{1.1}) enjoy  Lipschitz estimates  and are $C^{\infty}$ in each "strip" of $K$ up to the flat interface.
We present our analysis in the case $n\geq 3$ although our result can be extended easily to the case $n=2$. In fact the two dimensional case is in some sense easier to treat since in this case Dong and Kim in  \cite{DK} have proved existence, uniqueness and pointwise estimates of the Green's function in $\Omega$.

We want to point out that our result holds also if $K$ contains less regular  $C^{1,\alpha}$ interfaces (see Remark 5.1) but we think that the treatment  of this case would only require tedious and long technicalities and calculations.
On the other hand the flatness assumption of a portion of the interface  is not too restrictive since it includes for example a partition of
 $\Omega$ with polyhedral domains $D_j$ which appear in any numerical scheme used for the effective reconstruction of the admittivity.
Besides, this assumption   allows to derive H\"{o}lder quantitative estimates of unique continuation of solutions  to equation (\ref{1.1}) and
consequently  a better dependence  of the constant $C$ on $N$ .

The plan of the paper is the following: In Section 2 we introduce notation and the main assumptions and we state our main result (Theorem 2.1).
In Section 3 we collect all the results needed in order to prove Theorem 2.1.
 In Proposition 3.1 and Proposition 3.2 we state some known results concerning respectively the regularity of solutions of equation (\ref{1.1})
and the existence of the Green's function in the case of continuous admittivities. In the key Proposition 3.3 we prove the existence of
singular solutions in $K$ and we investigate their asymptotic  behaviour near the flat discontinuity interface.
  In Theorem 3.4 we show that the solutions of (\ref{1.1}), which due to the particular structure of $\gamma$ are piecewise analytic in $\Omega$,
can be extended analytically through the flat interfaces.  This property allows us in Proposition 3.5 to derive optimal quantitative estimates of unique continuation for solutions of equation (\ref{1.1}).  In Section 4 we give the proof of Theorem 2.2. In Section 5 we give some final remarks about generalizations of our result and, finally, the Appendix contains the statement of Caccioppoli inequality  (Proposition 6.1), the proof of Theorem 3.4, and the generalization of Alessandrini's identity to the complex case.

\section{Main result}\label{sec2}
\subsection{Notation and main assumptions}\label{subsec2.1}
For every $x\in\RR^n$ let us set $x=(x^\prime,x_n)$ where $x^\prime\in\RR^{n-1}$ for $n\geq 3$.
With $B_R(x)$ and $B_R^\prime(x^\prime)$ we will denote respectively the open ball in $\RR^n$ centered at $x$ of radius $R$
and the ball in $\RR^{n-1}$ centered at $x^\prime$ of radius $R$; $B_R(0)$ and $B_R^\prime(0)$ will be denoted by $B_R$  and $B^\prime_R$.
We will also use the following notations $\RR^n_+=\{(x^\prime,x_n)\in\RR^n\,:\,x_n>0\}$,
$\RR^n_-=\{(x^\prime,x_n)\in\RR^n\,:\,x_n<0\}$, $B^+_R=B_R\cap R^n_+$ and $B^-_R=B_R\cap R^n_-$.

We will denote by $D^\beta_x$ the derivative corresponding to a multiindex
$\beta=(\beta_1,\ldots,\beta_n)$ and by $D^{\beta^\prime}_{x^\prime}$
the partial derivative corresponding to the multiindex $\beta^\prime=(\beta_1,\dots,\beta_{n-1},0)$, while we will write
$\frac{\der}{\der x_h}$, for $h\in\{1,\ldots,n\}$, for the partial derivative with respect to $x_h$ and $\frac{\der}{\der \nu}$
the partial derivative in direction $\nu$.

Let $\om$ be a bounded domain in $\RR^n$. We shall say that 
$\der\om$ is Lipschitz continuous  with constants $r_0, L>0$ if $\forall P\in\der\om$ there exists a rigid transformation of coordinates such that $P=0$ and
\[
\om\cap B_{r_0}=\{(x^\prime,x_n)\in B_{r_0}\, :\, x_n> \phi(x^\prime)\},
\]
where $\phi$ is a Lipschitz continuous function on $B_{r_0}$ with $\phi(0)=0$ and
\[
\|\phi\|_{C^{0,1}\left(B_{r_0}^\prime\right)}\leq L r_0.
\]

Our main assumptions are:

\noindent{\bf (H1)}  $\om\subset\RR^n$ is a bounded domain such that
\[
|\om|\leq Ar_0^n,
\]
and
\[
\der\om \mbox{ is of Lipschitz class with constants } r_0, L.
\]
\noindent{\bf (H2)} The complex conductivity $\gamma$ satisfies
\begin{equation}\label{2.3}
\Re\gamma \geq \lambda^{-1},\quad |\gamma|\leq \lambda\quad\mbox{in}\quad\om
\end{equation}
for some $\lambda\geq 1$, and is of the form
\[
\gamma(x)=\sum_{j=1}^N\gamma_j \mathds{1}_{D_j}(x),
\]
where $\gamma_j$ are for $j=1,\ldots,N$ unknown complex numbers and $D_j$ are known open sets in $\RR^n$ which satisfy the following conditions

\noindent{\bf (H3)}
$D_j$, $j=1,\ldots,N$ are connected and pairwise nonoverlapping such that $\cup_{j=1}^N \ov{D}_j=\ov{\om}$,
$\der D_j$, $j=1,\ldots,N$ are of Lipschitz class with constants $r_0$, $L$.

We also assume that there exists one region, say $D_1$ such that $\der D_1\cap\der\om$ contains an open flat portion $\Sigma_1$.
For every $j\in\{2,\ldots,N\}$ there exist $j_1,\ldots,j_M\in\{1,\ldots,N\}$ such that
\[
D_{j_1}=D_1,\quad D_{j_M}=D_j
\]
and, for every $k=1,\ldots,M$
\[
\der D_{j_{k-1}}\cap \der D_{j_k}
\]
contains a flat portion $\Sigma_k$ such that
\[
\Sigma_k\subset\om,\quad \forall k=2,\ldots,M.
\]

Furthermore, there exists $P_k\in\Sigma_k$ and a rigid transformation of coordinates such that $P_k=0$ and
\begin{eqnarray*}
  \Sigma_k\cap B_{r_0/3} &=& \{x\in B_{r_0/3}\,:\, x_n=0\}, \\
  D_{j_k}\cap B_{r_0/3} &=& \{x\in B_{r_0/3}\,:\, x_n>0\}, \\
  D_{j_{k-1}}\cap B_{r_0/3} &=& \{x\in B_{r_0/3}\,:\, x_n<0\}.
\end{eqnarray*}

For simplicity we will call $D_{j_1},\ldots,D_{j_M}$ \textit{a chain of domains connecting $D_1$ to $D_j$}.
\vspace{0.1cm}

In the following we will introduce a number of constants that we will always denote by $C$. The values of this constants might
differ from one line to the other. We will write explicitly which a priori parameters each constant depends on.

Consider the problem
\begin{equation}\label{2.5}
    \left\{\begin{array}{rcl}
             \mbox{div}\left(\gamma\nabla u\right)&= & 0\quad\mbox{in}\quad\om \\
             u&= &f  \quad\mbox{on}\quad\der\om,\\
           \end{array}
    \right.
\end{equation}
where $\om$ and $\gamma$ satisfy assumptions {\bf (H1)-(H3)} and $f\in H^{1/2}(\der\om)$.
Observe that, by assumption (\ref{2.3}), applying Lax-Milgram Theorem, there exists a unique solution $u\in H^1(\om)$ of problem (\ref{2.5}).

For $f\in H^{1/2}(\der\om)$ define the Dirichlet-to-Neumann map $\Lambda_\gamma:H^{1/2}(\der\om)\rightarrow H^{-1/2}(\der\om)$ given by
\[
\Lambda_\gamma(f)=\gamma\frac{\der u}{\der \nu}_{|_{\der\om}}.
\]
Let us notice that $\Lambda_\gamma$ can be identified by the sesquilinear form
on $H^{1/2}(\der\om)\times H^{1/2}(\der\om)$ defined by
\[
<\Lambda_\gamma f, \psi>=\int_{\om}\gamma\nabla u\cdot \nabla \ov{v}, \quad \forall f,\psi\in H^{1/2}(\der\om)
\]
where $u$ is solution to problem (\ref{2.5}) and $v$ is any $H^1(\om)$ function such that $v_{|_{\der\om}}=\psi$.
%
\begin{teo}\label{teo2.1}
%
Let $\om$ satisfy assumption {\bf (H1)}. Let $\gamma^{(k)}$, $k=1,2$ be two complex piecewise constant functions of the form
\[
\gamma^{(k)}(x)=\sum_{j=1}^N\gamma^{(k)}_j \mathds{1}_{D_j}(x),
\]
where $\gamma^{(k)}$  satisfy for $k=1,2$ assumption {\bf (H2)} and $D_j, j=1,\ldots,N$ satisfy assumption
{\bf (H3)}.

Then there exists a positive constant $C=C(r_0,L,A,n,N,\lambda)$ such that
\begin{equation}\label{stabil}
\|\gamma^{(1)}-\gamma^{(2)}\|_{L^\infty(\om)}\leq C\|\Lambda_1-\Lambda_2\|_{\mathfrak{L}(H^{1/2}(\der\om),H^{-1/2}(\der\om))},
\end{equation}
where $\Lambda_i=\Lambda_{\gamma^{(k)}}$ for $k=1,2$.
\end{teo}
%
\section{Preliminary results}\label{sec3}
%
Let us now state some results that will be useful in the
proof of Theorem \ref{teo2.1}.

The first ingredient is a regularity estimate for solutions to the admittance equation in stratified media.
In order to get such a  regularity estimate we interpret equation $\mbox{div}(\gamma\nabla u)=0$, for a complex valued coefficient
$\gamma$, as a $2\times2$ differential system for real valued functions.
If we denote by $u^{(1)}=\Re u$ and $u^{(2)}=\Im u$ we have that the vector valued function $(u^{(1)},u^{(2)}):\om\rightarrow\RR^2$
satisfies the system
\begin{equation}\label{3.1}
\frac{\der}{\der x_h}\left(c_{lj}^{hk}\frac{\der u^{(j)}}{\der x_k}\right)=0,\quad l=1,2.
\end{equation}
where we used the convention of repeated index summation and where
\[c_{lj}^{hk}=\sigma\delta_{hk}\delta_{lj}-\ep\delta_{hk}(\delta_{l1}\delta_{j2}-\delta_{l2}\delta_{j1}),\]
for $l,j,h,k\in{1,2}$ and with $\sigma=\Re\gamma$ and $\ep=\Im\gamma$.
By assumption (\ref{2.3}), system (\ref{3.1}) satisfies the strong ellipticity condition
\[
    \lambda^{-1}|\xi|^2\leq c_{lj}^{hk} \xi_h^l\xi_k^j \leq \lambda |\xi|^2.
\]

For this  type of systems, Li and Nirenberg proved in \cite{LN} a regularity result that we state here for our particular equation.
The following proposition is a special case of Proposition 1.6 in \cite{LN}.

\begin{prop}\label{prop3.1}
Let $\gamma^{(1)}$ and $\gamma^{(2)}$ be two complex constants satisfying (\ref{2.3}). Let $r>0$ and let $h$, $g_1,\ldots,g_n$ be complex valued functions
of class $C^\infty\left(\ov{B}^\pm_r\right)$.

 Let $v\in H^1(B_r)$ be a solution to
\[
 \mbox{div}\left(\left(\gamma^{(1)}+(\gamma^{(2)}-\gamma^{(1)})\mathds{1}_{B^+_r}\right)\nabla v\right)=h+\mbox{div}\,g\quad\mbox{in}\quad  B_r,
 \]
 where $g=(g_1,\ldots,g_n)$.

Then, for every multi-index $\beta^\prime$, $D^{\beta^\prime}_{x^\prime} v\in C^0(B_r)$ and $v\in C^\infty(\overline{B}^\pm_r)$.
Moreover,  for every $\delta>0$ and  $k\geq0$,
\[
  \|v\|_{C^k(\overline{B}^\pm_{(1-\delta)r})}\!\!\leq C\left(\|v\|_{L^2\left(B_r\right)} +\sum_{|\beta|<\tilde{k}-1}\!\!\!\!\|D^\beta_{x^\prime}h\|_{L^2\left(B_r\right)}
  \vphantom{\|g\|_{C^\alpha\left(\overline{B}^+_{2\rho}(x)\right)}}+\sum_{|\beta|<\tilde{k}}\!\!\!\!\|D^\beta_{x^\prime}g\|_{L^2\left(B_r\right)}\right)
 \]
 where $\tilde{k}=k+ \left[\frac{n-1}{2}\right]+2$ and $C=C(\delta,k,n,\lambda)$.
\end{prop}
Note that, by $f\in C^\infty(\ov{B}^\pm_r)$ we intend that $f$ is separately  $C^\infty(\ov{B}^{\, +}_r)$ and
$C^\infty(\ov{B}^{\, -}_r)$.

Observe that, as a consequence of this result, since $v$ is continuous in $B_r$ and $\nabla v$ is bounded separately in $B_{(1-\delta)r}^+$
and in $B_{(1-\delta)r}^-$, then
$\nabla v\in L^\infty(B_{(1-\delta)r})$ and
\begin{equation}\label{stimagradiente}
    \|\nabla v\|_{L^\infty(B_{(1-\delta)r})}\!\!\leq C\left(\|v\|_{L^2\left(B_r\right)} +\!\!\!\sum_{|\beta|<\tilde{k}-1}\!\!\!\!\|D^\beta_{x^\prime}h\|_{L^2\left(B_r\right)}
  +\!\!\!\sum_{|\beta|<\tilde{k}}\!\!\|D^\beta_{x^\prime}g\|_{L^2\left(B_r\right)}\right)
\end{equation}
for $\tilde{k}=3+ \left[\frac{n-1}{2}\right]$.

These regularity estimates can be extended to $C^{1,\alpha}$ interfaces and less regular $h$ and $g$. (see Theorem 1.1 in \cite{LN}).

Our proof of Lipschitz stability estimates follows the approach used by Alessandrini and Vessella
 for the conductivity equation (\cite{AV}).
In their proof a crucial role is played by the Green function for the conductivity equation with bounded leading coefficient.
In our case, to our knowledge, the existence of a Green function in the whole domain $\om$ is not known for an $L^\infty$
complex coefficient in dimension $n\geq 3$. Existence of such a Green function is established in \cite{DM} for the $2$-dimensional
case  or for uniformly continuous coefficients in any dimension (see also \cite{HK} for a generalization of such a result to unbounded domains).

A Green matrix for a strongly elliptic operator
\[
Lu=-\frac{\der}{\der x_\alpha}(A^{\alpha\beta}_{lj}(x) \frac{\der u^j}{\der x_\beta})
\]
is a matrix valued function $\Gm\,:\, \left\{(x,y)\,:\, x,y\in\om,\, x\neq y\right\}\to\RR^{n\times n}$
such that
\begin{eqnarray*}
    L\Gm(\cdot,y)&=&\delta_y\,\mbox{Id}\quad\mbox{in}\,\om\\
\nonumber \Gm(\cdot,y)&=&0\quad\mbox{on}\,\der\om
\end{eqnarray*}
where $\delta_y$ is the Dirac distribution concentrated at $y$ and Id is the identity matrix in $\RR^n$,

If we have a Green matrix $\Gm=\{G_{ij}\}_{i,j=1}^n$ for system (\ref{3.1}), it is easy to see that the first row of $\Gm$, $G=(G_{11},G_{21})$ interpreted as a complex valued function, is a Green function for the operator $\mbox{div}(\gamma\nabla\cdot )$ in $\om$, in the sense that
\[
    \int_\om\gamma\nabla G(\cdot,y)\nabla\Phi=\Phi(y)
\]
for any complex valued function $\Phi\in C^\infty_0(\om)$.

Let us now state an existence result for the Green function for equation $\mbox{div}(\tilde{\gamma}\nabla u)=0$ for a continuous complex valued
coefficient $\tilde{\gamma}$.

\begin{prop}\label{prop3.2}
Let $\tilde{\gamma}\in C^0(\overline{\om})$ satisfy assumption (\ref{2.3}) and let $d_x=\mbox{dist}(x,\der\om)$. There exists a unique function
$\tilde{G}(x,y)$ continuous in $\left\{(x,y)\in\om\times\om\,:\, x\neq y\right\}$, locally integrable with respect to $y$ per every $x\in\om$
and such that, for every $f\in C^\infty_0(\om)$ the function
\[
u(x)=\int_\om \tilde{G}(x,y)f(y)\,dy
\]
belongs to $H^1(\om)$ and satisfies
\[-\mbox{div}\left(\tilde{\gamma}\nabla u\right)=f,\]
in the weak sense. Moreover,
\[
    \int_\om \tilde{\gamma}\nabla\tilde{G}(\cdot,y)\nabla\phi=\phi(y),\quad\mbox{for every }\phi\in C^\infty_0(\om),
\]
and, for every $\eta\in C^\infty_0(\om)$ such that $\eta\equiv 1$ in $B_r(y)$ for some $r<d_y$,
 \[
    (1-\eta)\tilde{G}(\cdot,y)\in H^1(\om).
 \]
 Furthermore
 \begin{equation}\label{3.15}
    \tilde{G}(x,y)=\tilde{G}(y,x)\quad\mbox{for every }x,y\in\om
 \end{equation}
 and
\[
    \left\|\tilde{G}(\cdot,y)\right\|_{H^1(\om\setminus B_r(y))}\leq C r^{1-\frac{n}{2}},\quad\mbox{for every }r<d_y/2.
 \]
\end{prop}

\textit{Proof.} This result follows from Theorem 5.4 in \cite{HK} and the observation that, given a Green matrix for system (\ref{3.1}),
it is possible to get existence of the Green function for equation $\mbox{div}\left(\tilde{\gamma}\nabla u\right)=0$.
The symmetry result (\ref{3.15}) follows from 
Theorem 1 in \cite{DM}.
\qed

These results on the Green function have not been extended to $L^\infty$ coefficients, hence we cannot use a Green function for our problem.
For this reason we will construct some solution of our equation that has the same behavior of a Green function, but only for $y$ in a
certain special subset of $\om$.
 Before doing this we need to extend our original domain $\om$ to a $\om_0$ by adding an open cylinder $D_0$ whose basis is the flat
portion $\Sigma_1$ of $\der\om\cap\overline{D}_1$ and with height greater than $r_0$.
Let $K_0=\{x\in D_0\,:\, \mbox{dist}(x,\Sigma_1)\geq r_0/2\}$.
If we set $\om_0=\om\cup D_0$ then $\der\om_0$ is Lipschitz continuous.

We extend any complex coefficient $\gamma$ defined in $\om$ by setting it equal to $1$ in $D_0$. For simplicity we will still denote this extension with $\gamma$.


Let us consider any subdomain in $\om$ and let us consider the chain of  domains connecting it to $D_1$ (see assumption \textbf{(H3)}).
For simplicity let us rearrange the indices of subdomains so that this chain corresponds to
$D_0, D_1,\ldots, D_M$, $M\leq N$. Let us denote by $\mathcal{S}=\cup_{j=0}^M \overline{D}_j$ and consider
a connected subset $K\subset\mathcal{S}$ with Lipschitz boundary such that
$\overline{K}\cap \der D_j=\Sigma_j\cup\Sigma_{j+1}$ for $j=1,\ldots,M$, $K\supset K_0$ and
 $\mbox{dist}(K,\der \mathcal{S}\setminus \{\Sigma_{M+1}\cup D_0\})>r_0/16$ .

\vspace{1em}

Let us denote by $\Gamma(x,y)$ the standard fundamental solution for the Laplace equation, given by
\[
\Gamma(x,y)=\frac{1}{(n-2)\omega_n}|x-y|^{2-n}\quad\mbox{in }\RR^n, n\geq 3.
\]
where $\omega_n$ is the volume of the unit ball in $n$ dimensions.

Let $\gamma,\delta\in\CC$. Then a straightforward calculation shows that
\[
    \Gamma_{\gamma,\delta}(x,y)=\left\{\begin{array}{ccl}
    \frac{1}{\gamma}\Gamma(x,y)+s\Gamma(x,y^*)&\mbox{if}& x_n>0, y_n>0,\\
    \left(\frac{1}{\gamma}+s\right)\Gamma(x,y)&\mbox{if}&x_n\cdot y_n<0,\\
    \frac{1}{\delta}\Gamma(x,y)+t\Gamma(x,y^*)&\mbox{if}&x_n<0, y_n<0,
    \end{array}\right.
\]
where $y^*=(y_1,\ldots,y_{n-1},-y_n)$, $s=\frac{\gamma-\delta}{\gamma(\gamma+\delta)}$ and $t=\frac{\delta-\gamma}{\gamma(\gamma+\delta)}$,
is a fundamental solution for the differential operator
\begin{equation}\label{3.17}
\mbox{div}\left(\left(\delta \mathds{1}_{\RR^n_-}+\gamma\mathds{1}_{\RR^n_+}\right)\nabla\cdot\right)\quad\mbox{in}\quad \RR^n.
\end{equation}

\begin{prop}\label{prop3.3}
Let $\gamma$ satisfy assumptions {\bf (H1)-(H3)} in $\om_0$ and let $\mathcal{S}$, $K_0$ and $K$ be defined as above.
For $y\in K$ there exists a unique function $G(\cdot,y)$, continuous in $\om\setminus\{y\}$
such that
\[
\int_\om\gamma\nabla G(\cdot,y)\cdot\nabla \phi=\phi(y),\quad \forall\phi\in C^\infty_0(\om).
\]

Furthermore,
\begin{equation}\label{3.20}
    \|G(\cdot,y)\|_{H^1(\om\setminus B_r(y))}\leq C r^{1-n/2},\quad\forall r<d_y/2,
\end{equation}
and
\begin{equation}\label{3.sym}
G(x,y)=G(y,x)\quad\mbox{for every }x,y\in K.
\end{equation}

Let $D_l$ and $D_{l+1}$ be two subdomains of $\mathcal{S}$ such that $\der D_{l+1}\cap\der D_l$ contains a flat portion $\Sigma_{l+1}$ satisfying assumption (H3). Let us fix the origin at $P_{l+1}\in \Sigma_{l+1}$ and let $\nu$ be the outer normal to $D_l$ at the origin. Let $\overline{y}=-r\nu$ for some $r\in (0,r_0/6)$ and let $\overline{x}\in B_{r_0/6}\cap D_{l+1}$. There exists a constant $C=C(r_0,\lambda,n,A,L)$ such that
\begin{eqnarray}\label{3.21}
   &\left|G(\overline{x},\overline{y})-\frac{2}{\gamma_{l}+\gamma_{l+1}}\Gamma(\overline{x},\overline{y})\right| \leq C &\nonumber\\
   &&  \\
&\left|\nabla_x G(\overline{x},\overline{y})-\frac{2}{\gamma_{l}+\gamma_{l+1}}\nabla_x\Gamma(\overline{x},\overline{y})\right| \leq C & \nonumber
\end{eqnarray}
\end{prop}

\textit{Proof of Proposition \ref{prop3.3}.} Let us construct $G$ by taking advantage of the fundamental solution for operator (\ref{3.17}).
If $y\in K$, there is a couple of contiguous domains of $\mathcal{S}$ such that $y\in \mathcal{S}\cap (\overline{D}_l\cup\overline{D}_{l+1})$ and
$\mbox{dist}(y,D_j)\geq \frac{r_0}{6}$ for every $j\in\{0,\ldots,M\}\setminus \{l, l+1\}$. Let us fix the origin at the point
$P_{l+1}\in\Sigma_{l+1}$. Let us denote by $\Gamma_l(x,y)=\Gamma_{\gamma_l,\gamma_{l+1}}(x,y)$ and by
$\tilde{\gamma}=\gamma-\gamma_l\mathds{1}_{\RR^n_-}-\gamma_{l+1}\mathds{1}_{\RR^n_+}$
and let
\[G(x,y)=\Gamma_l(x,y)+w(x,y),\]
where $w$ is solution to
\[
   \left\{\begin{array}{rcl}
             \mbox{div}\left(\gamma \nabla_x w(\cdot,y)\right)&= &\mbox{div}\, h\quad\mbox{in}\quad\om \\
             w(\cdot,y)&= &-\Gamma_l(\cdot,y)  \quad\mbox{on}\quad\der\om,\\
           \end{array}
    \right.
 \]
 where $h=\tilde{\gamma}\nabla_x\Gamma_l(\cdot,y)$.
This problem has a unique solution because $-\Gamma_l(\cdot,y)\in H^{1/2}(\der\om)$ and since $\tilde{\gamma}=0$ in $D_l\cup D_{l+1}$,
$\mbox{div}\,h\in H^{-1}(\om)$.

Now, in order to get (\ref{3.21}) consider
\[
w_0(x)=w(x)+\tilde{\Gamma}_l(x,y),
\]
where $\tilde{\Gamma}_l(x,y)=\phi(x)\Gamma_l(x,y)$ with $\phi\in C^\infty(\overline{\om})$, $\phi=0$
in $B_{r_0/3}(y)$ and $\phi=1$ in $\om\setminus\overline{B_{2r_0/3}(y)}$.

The function $w_0$ satisfies
 \begin{equation}\label{3.22}
   \left\{\begin{array}{rcl}
             \mbox{div}\left(\gamma \nabla w_0\right)&= &\mbox{div}(h_0+h)\quad\mbox{in}\quad\om \\
             w_0&= &0  \quad\mbox{on}\quad\der\om,\\
           \end{array}
    \right.
 \end{equation}
 where $h_0=\gamma \nabla_x \tilde{\Gamma}_l(\cdot,y)$.
Multiplying equation (\ref{3.22}) by $\overline{w}_0$, integrating by parts, using Schwartz inequality and the fact that $\tilde{\gamma}=0$ in $D_l\cup D_{l+1}$ and $\nabla \tilde{\Gamma}_l=0$ in $B_{r_0/3}(y)$ we get
 \begin{eqnarray*}
   \int_\om|\nabla w_0(\cdot,y)|^2 \,dx &\leq& C\left|\int_\om \tilde{\gamma}\nabla \Gamma_l(\cdot, y)\cdot \nabla \overline{w}_0\right| \\
    &+&C \left|\int_\om \gamma\nabla \tilde{\Gamma}_l(\cdot, y)\cdot \nabla \overline{w}_0\right|\\
    &\leq& C\|\nabla w_0\|_{L^2(\om)}^{1/2},
 \end{eqnarray*}
 where $C=C(r_0,n,\lambda,A,L)$.
By Poincar\'{e} inequality this also implies that
\begin{equation}\label{E3}
    \|w_0\|_{H^1(\om)}\leq C
\end{equation}
Estimate (\ref{3.20}) follows then immediately from  (\ref{E3}) and from the behavior of $\Gamma_l$.

Since $h=0$ in $D_l\cup D_{l+1}$ and $h_0$ belong to $C^\infty\left(B_{\frac{r_0}{3}}^\pm\right)$, we can apply to function $w_0$ estimate (\ref{stimagradiente}) and get
 \begin{eqnarray}\label{*}
   \|\nabla w_0\|_{L^\infty(B_{r_0/6})} &\leq& C\left(\|w_0\|_{L^2(B_{r_0/3})}+\sum_{|\beta|\leq \tilde{k}}\|D^\beta_{x^\prime} h_0\|_{L^2(B_{r_0/3})}\right)\\
   ,\nonumber
 \end{eqnarray}
 where  $\tilde{k}=3+\left[\frac{n-1}{2}\right]$.

By the previous inequality, (\ref{E3}) and (\ref{*}) we get that
\[
\|\nabla w_0\|_{L^\infty(B_{r_0/6})}\leq C\quad\mbox{and}\quad \|w_0\|_{L^\infty(B_{r_0/6})}\leq C
\]
Finally, since $w=w_0-\tilde{\Gamma}_l$ and $\|\nabla\tilde{\Gamma}_l\|_ {L^\infty(B_{r_0/6})},\|\tilde{\Gamma}_l\|_ {L^\infty(B_{r_0/6})}\leq C$,  (\ref{3.21}) follows.

Symmetry (\ref{3.sym}) of $G$ in $K$ follows by standard arguments based on integration by parts
(see  for example \cite[Theorem 13, p. 35]{E} ).\qed

Solutions to our equation are harmonic in each of the subdomains $D_j$, hence they are piecewise analytic in $\om$. Now we want to show that is that we can analytically extend each analytic portion through the flat interface. This property of $u$ will allow us in Proposition 3.5 to derive Hoelder estimates of unique continuation in $K$.
 More precisely we will prove the following

\begin{teo}\label{teo3.3}
Let $D_1,\ldots D_M$ be a chain and let $u$ be a solution to
\[
    \mbox{div}\left(\gamma\nabla u\right)=0\quad\mbox{in}\quad D_j\cup D_{j+1}
\]
for some $j\in{1,2,\ldots,M}$. Then there exist two positive constants $C_1=C_1(\lambda,n)$ and $C=C(\lambda,r_0,A,L,n)$
 such that $u_{|_{D_j}}$ can be extended
by a function $\tilde{u}$ analytic in the set $D_j\cup E_{j+1}^{(C_1)}$, where
\[E_{j+1}^{(C_1)}=\left\{x\in D_{j+1}\,:\, \mbox{dist}(x,B_{\frac{r_0}{2}}(x_0)\cap\Sigma_{j+1})<\frac{r_0}{4 C_1}\right\},\]
and
\begin{equation}\label{3.32}
    \|\tilde{u}\|_{L^\infty\left(D_j\cup E_{j+1}^{(C_1)}\right)}\leq C\|u\|_{L^2(D_j\cup E_{j+1}^{(2)})}.
    \end{equation}
\end{teo}
The proof of this Theorem is contained in the Appendix.

%
%
\begin{prop}\label{prop3.4}
Let  $K$ and $K_0$ as before, and let $v\in H^1(K)$
be a solution to
\[\mbox{div}\left(\gamma\nabla v\right)=0\quad\mbox{in}\quad K,\]
such that
\[
    \|v\|_{L^\infty(K_0)}\leq C\ep_0r_0^{2-n},
\]
and
\[
    |v(x)|\leq C(\ep_0+E_0)r_0^{1-\frac{n}{2}}\mbox{dist}(x,\Sigma_{M+1})^{1-\frac{n}{2}}
    \quad\mbox{for every}\quad x\in K.
\]
Then
\begin{equation}\label{3.35}
    |v(\tilde{x})|\leq C \left(\frac{\ep_0}{E_0+\ep_0}\right)^{\tau^{(M+1)N_1}\delta_1^{M+1}\tau_r}
    (E_0+\ep_0)r_0^{2-n}\left(\frac{r}{r_0}\right)^{\left(1-\frac{n}{2}\right)(1-\tau_r)},
\end{equation}
where $\tilde{x}=P_{M+1}-2r\nu\left(P_{M+1}\right)$, $r\in\left(0,\frac{3}{40}r_0\right)$,
 $\tau=\frac{\ln(4/3)}{\ln 4}$, $\delta_1\in(0,1)$, $N_1$ and $C$  depend on $r_0,L,A,n,\lambda$,
 and $\tau_r=\frac{\ln\left(\frac{3r_1-r}{3r_1-2r}\right)}{\ln\left(\frac{3r_1-r}{r_1}\right)}$,
for $r_1=r_0/4C_1$ and $C_1$ as in Theorem \ref{teo3.3}.
\end{prop}

\textit{Proof.}  By Theorem \ref{teo3.3}, the function $v_{|_{D_0}}$ can be extended analytically to a function $v_0$ on $D_0\cup E_1^{(C_1)}$
such that ${v_0}_{|_{D_0}}=v_{|_{D_0}}$ and
\[\|v_0\|_{L^\infty(D_0\cup E_1^{(C_1)})}\leq C (E_0+\ep_0)r_0^{2-n}.\]
Let $r_1=\frac{r_0}{4 C_1}$ where $C_1$ and $E_1^{(C_1)}$ are the same as in Theorem \ref{teo3.3}.
Note  that, $C_1>>4$ because $\lambda>1$.

Let us consider the sphere $B_{4r_1}$ of radius $4r_1=\frac{r_0}{C_1}\leq\frac{r_0}{4}$ strictly contained in $K_0$.
Let $B_{3r_1}$ and $B_{r_1}$ be spheres concentric to $B_{4r_1}$ and of radius $3r_1$ and $r_1$ respectively.
Let $P\in \Sigma_1$ such that $\mbox{dist}(P,P_1)<\frac{r_0}{2}$. Let us construct a chain of spheres of radius $r_1$ such that
the first is $B_{r_1}$, all the spheres are externally tangent and the last one is centered at $P-2r_1\nu_1$ where $\nu_1$ is the
exterior normal vector to $\Sigma_1$. We choose this chain so that the spheres of radius $4r_1$ concentric with those of the
chain are contained in $D_0\cup E_1^\lambda$. Such a chain has a finite number of spheres that is certainly smaller than $N_1=\frac{|\om|}{c_nr_1^n}+1$.

By the three sphere inequality (see, for example \cite{ADB}) we have that
\[
\|v_0\|_{L^\infty\left(B_{3r_1}(P-2r_1\nu_1)\right)}\leq C\ep_0^{\tau^{N_1}}(E_0+\ep_0)^{1-\tau^{N_1}}r_0^{2-n},
\]
where $C=C(\lambda,r_0,n)$.

In particular, for every $P\in\Sigma_1$ such that $\mbox{dist}(P,P_1)\leq \frac{r_0}{2}$, we have
\[
\|v_0\|_{L^\infty(B_{r_1}(P))}\leq C\ep_0^{\tau^{N_1}}(E_0+\ep_0)^{1-\tau^{N_1}}r_0^{2-n}.
\]
Since $\Re v_0$ and $\Im v_0$ are harmonic,
\[
\|\nabla v_0\|_{L^\infty\left(B_{\frac{r_1}{2}}(P)\right)}\leq Cr_1^{-1}\ep_0^{\tau^{N_1}}(E_0+\ep_0)^{1-\tau^{N_1}}r_0^{2-n}.
\]
This also implies that the Cauchy data of $v$ on $\Sigma_1$ are small,
\begin{eqnarray*}
  \|v\|_{L^\infty\left(\Sigma_1^-\cap B^\prime_{\frac{r_0}{2}}(P_1)\right)} &\leq& C\ep_0^{\tau^{N_1}}(E_0+\ep_0)^{1-\tau^{N_1}}r_0^{2-n}, \\
  \|\nabla v\|_{L^\infty\left(\Sigma_1^-\cap B^\prime_{\frac{r_0}{2}}(P_1)\right)} &\leq& C\ep_0^{\tau^{N_1}}(E_0+\ep_0)^{1-\tau^{N_1}}r_0^{1-n}.
\end{eqnarray*}
Let us now consider $v_{|_{D_1}}$. Due to the transmission conditions
\begin{eqnarray*}
  \|v\|_{L^\infty\left(\Sigma_1^{+}\cap B^\prime_{\frac{r_0}{2}}(P_1)\right)} &\leq& C\ep_0^{\tau^{N_1}}(E_0+\ep_0)^{1-\tau^{N_1}}r_0^{2-n}, \\
  \left\|\frac{\der v}{\der \nu_1}\right\|_{L^\infty\left(\Sigma_1^+\cap B^\prime_{\frac{r_0}{2}}(P_1)\right)} &\leq&
C\ep_0^{\tau^{N_1}}(E_0+\ep_0)^{1-\tau^{N_1}}r_0^{1-n}.
\end{eqnarray*}
By Trytten Cauchy estimates for solutions to elliptic equations (cfr. \cite{T} and \cite{ABRV})  we get that
\begin{equation}\label{trytten}
\int_{D_1\cap B_{\frac{r_0}{2}-\frac{r_0}{16}}} |\nabla v|^2 dx\leq C\ep_0^{2\tau^{N_1}\delta_1}(E_0+\ep_0)^{2(1-\tau^{N_1}\delta_1)}r_0^{2(2-n)}
\end{equation}
for some $\delta_1\in (0,1)$ depending only on the a priori data.

By standard regularity estimates for harmonic functions we get
\[
\|v\|_{L^\infty\left(D_1\cap B_{\frac{r_0}{2}-\frac{r_0}{8}}(P_1)\right)}\leq C\ep_0^{\tau^{N_1}\delta_1}(E_0+\ep_0)^{1-\tau^{N_1}\delta_1}r_0^{2-n}.
\]
Again we can apply the three sphere inequality considering $r_1\leq\frac{r_0}{16}$ with $B_{4r_1}\subset D_1$ and
$B_{r_1}\subset D_1\cap B_{\frac{r_0}{2}-\frac{r_0}{8}}$ and applying again Theorem \ref{teo3.3}, we get that the analytic extension
$v_1$ of $v_{|_{D_1}}$ satisfies
\begin{eqnarray*}
  \|v_1\|_{L^\infty\left(B_{r_1}(P)\right)} &\leq& C\ep_0^{\tau^{2{N_1}}\delta_1}(E_0+\ep_0)^{1-\tau^{2{N_1}}\delta_1}r_0^{2-n}, \\
  \|\nabla v_1\|_{L^\infty\left(B_{r_1}(P)\right)} &\leq& C\ep_0^{\tau^{2N_1}\delta_1}(E_0+\ep_0)^{1-\tau^{2N_1}\delta_1}r_0^{1-n}.
\end{eqnarray*}
for every $P\in\Sigma_2$ such that $\mbox{dist}(P,P_2)<\frac{r_0}{2}$.

Hence again by smallness of Cauchy data, transmission conditions and Trytten inequality (\ref{trytten}),
we get
\[
\left|v(x)\right|\leq C \ep_0^{\tau^{(M+1)N_1}\delta_1^{M+1}}(\ep_0+E_0)^{1-\tau^{MN_1}\delta_1^{M}}r_0^{2-n}.
\]
for every $x$ such that $x=P-2tr_1\nu(P_{M+1})$ where $P\in\Sigma_{M+1}$, $\mbox{dist}(P,P_{M+1})\leq \frac{r_0}{2}$ and $1<t<2$.

Now, let $\tilde{x}=P_{M+1}-r\nu(P_{M+1})$ where $r<r_1$. We can again use a three spheres inequality for the spheres centered at
$P_{M+1}-3r_1\nu(P_{M+1})$ and of radii $r_1$, $3r_1-2r$ and $3r_1-r$ and get (\ref{3.35}).\qed

\section{Proof of Theorem \ref{teo2.1}}\label{sec4}
Let $D_0, D_1,\ldots,D_M$ be the chain of domains such that
\[
\|\gamma^{(1)}-\gamma^{(2)}\|_{L^\infty(D_M)}=E:=\|\gamma^{(1)}-\gamma^{(2)}\|_{L^\infty(\om)}.
\]
Consider $\mathcal{S}$, $K$ and $K_0$ as defined in the previous section.

For $y\in K$, let $G_1(x,y)$ and $G_2(x,y)$ be the singular function related to $\gamma^{(1)}$ and $\gamma^{(2)}$, respectively,
whose existence and behavior  has been shown in Proposition \ref{prop3.3}.

Let $U_0=\om$ and $U_k$, for $k=1,\ldots,M$, be given by $U_k=\om\setminus \cup_{j=1}^k D_j$. Let also denote by $W_k=\cup_{j=0}^k D_j$.
Define, for $y,z\in K$,
\[
S_k(y,z):=\int_{U_k}(\gamma^{(1)}(x)-\gamma^{(2)}(x))\nabla_x G_1(x,y)\cdot\nabla_x G_2(x,z)\,dx
\]
Observe that, by (\ref{3.20}),
\[
    \left|S_k(y,z)\right|\leq C E\left(d(y)d(z)\right)^{1-\frac{n}{2}}, \quad\mbox{for every } y,z\in K,
\]
where $C$ depends on the a priori assumptions and $d(y)=\mbox{dist}(y, U_k)$.

The main point of the proof consists in showing that
\begin{prop}\label{prop4.1}
For every $y,z\in K$ we have $S_k(\cdot,z), S_k(y,\cdot)\in H^1_{loc}(K)$ and
\begin{equation}\label{4.3}
    \mbox{div}\left(\gamma^{(1)}\nabla S_k(\cdot,z)\right)=0,\quad \mbox{div}\left(\gamma^{(2)}\nabla S_k(y,\cdot)\right)=0\quad\mbox{in }K.
\end{equation}
\end{prop}
For the proof of Proposition \ref{prop4.1} we need the following  approximation result that we will prove later on.
\begin{lm}\label{lemma4.2}
Let $\gamma$ satisfy assumptions \textbf{(H1)-(H3)} and let $D_0$, $\om_0$, $U_k$ and $K$ as above, and,  for $y\in K$,
let $G(x,y)$ be the function defined in Proposition \ref{prop3.3}.
Assume that $\left\{\gamma_h\right\}_{h\in \NN}$ is a sequence of complex valued continuous coefficients that converges to $\gamma$ in $L^s(\om_0)$ for every
$s\in[1,+\infty)$ and let $G_h(x,y)$ be the Green's function for $\om_0$.
Then
\begin{equation}\label{4.4}
    \lim_{h\to+\infty}\sup_{y\in \tilde{K}}\|G_h(\cdot,y)-G(\cdot,y)\|_{H^1(U_k)}=0,
\end{equation}
for every $\tilde{K}\subset\subset K$.
\end{lm}

\textit{Proof of Proposition \ref{prop4.1}.}
We consider, for $j=1,2$ a regularization of $\gamma^{(j)}$ obtained by convolution with mollifiers
$\alpha_h\in C^\infty_0(\RR^n)$ such that $0\leq \alpha_h\leq 1$, $\alpha_h=0$ for $|x|>\frac{1}{h}$,
$\int_{\RR^n}\alpha_h=0$.
Then consider
\[
\gamma^{(j)}_h(x)=\int_{\RR^n}\gamma^{(j)}(y)\alpha_h(x-y)\,dy.
\]
Clearly $\gamma^{(j)}_h$ is smooth and $\Re \gamma^{(j)}_h\geq\frac{1}{\lambda}$ and $|\gamma^{(j)}_h|\leq \lambda$.

Let $G_{j,h}$ be the Green's function of $\mbox{div}\left(\gamma^{(j)}_h \nabla\cdot\right)$.

Let $\tilde{K}\subset\subset K$. By (\ref{4.4})
\begin{equation}\label{6.1}
    \sup_{\tilde{K}}\left\|G_{1,h}(\cdot,y)-G_1(\cdot,y)\right\|_{H^1(U_k)}\to 0 \quad\mbox{as }h\to+\infty.
\end{equation}
Set
\[S_{k,h}(y,z)=\int_{U_k}(\gamma^{(1)}-\gamma^{(2)})\nabla_x G_{1,h}(x,y)\nabla_x G_2(x,z)\,dx,\]
for $y,z\in K$.
Let $z\in K$; by using the symmetry of $G_{1,h}$ and differentiating under the integral sign, we have
\[
\mbox{div}\left(\gamma^1_h \nabla S_{k,h}(\cdot,z)\right)=0 \quad\mbox{in }K, \quad\mbox{for every }h.
\]
By (\ref{6.1}),
\begin{equation}\label{6.3}
    S_{k,h}(\cdot,z)\to S_k(\cdot,z)\quad\mbox{in }L^\infty(\tilde{K}).
\end{equation}
Let $K_1$ be such that $\tilde{K}\subset \mathring{ K_1}\subset \overline{K}_1\subset K$. By Caccioppoli inequality \ref{Caccio} we have that
\[\|S_{k,h}(\cdot,z)\|_{H^1(\tilde{K})}\leq C \|S_{k,h}(\cdot,z)\|_{L^2(K_1)}.\]
On the other hand, it is easy to check that
\[\|S_{k,h}(\cdot,z)\|_{L^2(K_1)}\leq C \left(\mbox{dist}(K_1,K)d(z)\right)^{1-\frac{n}{2}}.\]
This implies that, by considering a subsequence, $S_{k,h}(\cdot,z)$ converges weakly in $H^1(\tilde{K})$ to $S_k(\cdot,z)$ (by (\ref{6.3})).
In particular this implies that
\[
\mbox{div}\left(\gamma^{(1)}\nabla S_k(\cdot,z)\right)=0\quad\mbox{in }K_1.
\]
Since $K_1$ is any domain compactly contained in $K$, and since we can proceed in the same way with respect to $z$, we can conclude that
(\ref{4.3}) holds.

\textit{Proof of Lemma \ref{lemma4.2}}
Since $\gamma_h\to\gamma$ in $L^s(\om_0)$ for every $s\geq 1$, there exists a subsequence, that we continue to denote by $\{\gamma_h\}$ that converges a.e. to
$\gamma$ in $\om_0$.

For some $f\in C^\infty(\om_0)$, let $u_0$ and $u_h$ in $H^1_0(\om_0)$ be solutions in $\om$ to $\mbox{div}\left(\gamma \nabla u_0\right)=-f$ and $\mbox{div}\left(\gamma_h \nabla u_h\right)=-f$
respectively.
Observe that
\[\mbox{div}\left(\gamma_h\nabla(u_h-u_0)\right)=\mbox{div}\left((\gamma-\gamma_h)\nabla u_0\right)\quad\mbox{in }\om_0.\]
Multiplying the above equation by $\overline{u}_h-\overline{u}_0$ and integrating by parts we get
\[
\int_\om\gamma_h|\nabla(u_h-u_0)|^2=\int_\om (\gamma-\gamma_h)\nabla u_0\cdot\nabla (\overline{u_h-u_0}).
\]
Using the strong ellipticity condition and H\"{o}lder's inequality we get
\[
\int_{\om_0}|\nabla (u_h-u_0)|^2\leq \lambda\left(\int_{\om_0} |\gamma-\gamma_h|^2|\nabla u_0|^2\right)^{\frac{1}{2}}\left(|\nabla (u_h-u_0)|^2\right)^{\frac{1}{2}},
\]
hence
\[
\left\|\nabla (u_h-u_0)\right\|_{L^2(\om_0)}\leq \lambda \left(\int_{\om_0}|\gamma-\gamma_h|^2|\nabla u_0|^2\right)^{\frac{1}{2}}.
\]
By the dominated convergence theorem, we get
\[
\|\nabla(u_h-u_0)\|_{L^2(\om_0)}\to 0\mbox{ as }h\to 0
\]
and, by Poincar\'{e} inequality
\[
\|u_h-u_0\|_{H^1(\om_0)}\to 0 \mbox{ as }h\to 0.\]
Multiplying equation  $\mbox{div}\left(\gamma_h \nabla u_h\right)=-f$ by $\overline{u_h}$ we get
\[
\int_{\om_0}\gamma_h|\nabla u_h|^2=\int_{\om_0}f u_h.
\]
By strong ellipticity and Schwartz inequality we get
\[
\int_{\om_0}|\nabla u_h|^2\leq \lambda \|f\|_{L^2(\om_0)}\|u_h\|_{L^2(\om_0)}
\]
and, by Poincar\'{e} inequality
\[
\int_{\om_0}|\nabla u_h|^2\leq C\|f\|_{L^2(\om_0)}\|\nabla u_h\|_{L^2(\om_0)}.
\]
Hence
\[
\|u_h\|_{H^1(\om_0)}\leq C \|f\|_{L^\infty(\om_0)}\quad\mbox{for }h\geq 0.
\]

Functions $u_0$ and $u_h$ satisfy in $K$ the assumptions of Proposition \ref{prop3.1}, hence, for $h\geq 0$,
\[
\sup_{y\in \tilde{K}} |u_h|, |\nabla u_h| \leq C\left(\|u_h\|_{L^2(\om_0)}+\|f\|_{C^{2+\left[\frac{n-1}{2}\right]}(\om_0)}\right)\leq C\|f\|_{C^{2+\left[\frac{n-1}{2}\right]}(\om_0)}
\]
This implies that, for some subsequence of $\{u_h\}$
\begin{equation}\label{4.6}
    u_h\to u_0\quad\mbox{uniformly in }K
\end{equation}
and $u_0$ is continuous in $\tilde{K}$.

By the properties of functions $G$ and $G_h$ we have that
\[
u_0(y)=\int_{\om_0}G(z,y)f(z)\,dz\quad\mbox{for }y\in K,
\]
and
\[
u_h(y)=\int_{\om_0}G_h(z,y)f(z)\,dz\quad\mbox{for }y\in \om_0.
\]
By (\ref{4.6}), uniformly with respect to $y\in \tilde{K}$,
\[
    \int_{\om_0}G_h(z,y)f(z)\,dz\to \int_{\om_0}G(z,y)f(z)\,dz
\]
Let $\tilde{Q}$ be a smooth domain such that $U_k\subset\subset \tilde{Q}\subset\om_0$ such that
$\mbox{dist}(\tilde{Q},\tilde{K})>0$.

Now, let $x\in \tilde{Q}$; by symmetry of $G_h$ and $G$,
\begin{eqnarray*}
 \mbox{div}\left(\gamma_h\nabla G_h(x,\cdot)\right) &=0& \mbox{in }K,  \\
 \mbox{div}\left(\gamma\nabla G(x,\cdot)\right) &=0& \mbox{in }K.
\end{eqnarray*}
Again by Proposition \ref{prop3.1}, $G$ and $G_h$ satisfy a $C^{0,1}(\tilde{K})$ bound
uniformly with respect to $x\in \tilde{Q}$ and $h\in\NN$.

Hence
\[
\|G_h(\cdot,y)-G(\cdot,y)\|_{L^2(\tilde{Q})}
\]
satisfies a $C^{0,1}(\tilde{K})$ bound uniformly with respect to $h\in\NN$.

Thus, there exists a sequence $\{y_h\}\subset \overline{\tilde{K}}$ such that
\begin{equation}\label{s*}
\sup_{y\in \tilde{K}}\|G_h(\cdot,y)-G(\cdot,y)\|_{L^2(\tilde{Q})}=\|G_h(\cdot,y_h)-G(\cdot,y_h)\|_{L^2(\tilde{Q})}
\end{equation}
and $y_h\to\overline{y}\in \overline{\tilde{K}}$ (as a matter of fact this holds for some subsequence to which we restrict). Therefore
\begin{equation}\label{s**}
\|G_h(\cdot,y_h)-G_h(\cdot,\overline{y})\|_{L^2(\tilde{Q})}\to 0\quad\mbox{as }h\to+\infty.
\end{equation}
By (\ref{s*}) and (\ref{s**}) we have that, for every $f\in C^\infty(\om_0)$ such that $\mbox{supp}(f)\subset \tilde{Q}$, we have
\begin{eqnarray}\label{s*1}
 \nonumber 
 && \left|\int_{\om_0}\left(G_h(z,y_h)-G(z,\overline{y})\right)f(z)\,dz\right| \leq   \left|\int_{\om_0}\left(G_h(z,y_h)-G_h(z,\overline{y})\right)f(z)\,dz\right|  \\
 && \hphantom{aaaaaaaa} +\left|\int_{\om_0}\left(G_h(z,\overline{y})-G(z,\overline{y})\right)f(z)\,dz\right|\to 0\quad\mbox{ as }h\to+\infty.
\end{eqnarray}
 As a solution of $\mbox{div}\left(\gamma_h\nabla G_h(\cdot, y_h)\right)=0$ in $\om_0\setminus \overline{\tilde{K}}$, by \cite[Theorem 9.1, p.95]{C}
\begin{equation}\label{s*2}
 \|G_h(\cdot,y_h)\|_{W^{1,p}(\tilde{Q})}\leq C
 \end{equation}
for some $p>2$ and $C$ independent of $h$. Hence by Sobolev imbedding theorem, taking into account (\ref{s*1}) and (\ref{s*2}), up to subsequences,
\[
G_h(\cdot, y_h)\to G(\cdot, \overline{y})\quad\mbox{in }L^2(\tilde{Q}).
\]
%

Now observe that, for $h\to+\infty$
\begin{eqnarray*}
  \sup_{y\in\tilde{K}}\|G_h(\cdot,y)-G(\cdot,y)\|_{L^2(\tilde{Q})} = \|G_h(\cdot,y_h)-G(\cdot,y_h)\|_{L^2(\tilde{Q})}\leq\\
  \leq \|G_h(\cdot,y)-G(\cdot,\overline{y})\|_{L^2(\tilde{Q})}+\|G(\cdot,\overline{y})-G(\cdot,y_h)\|_{L^2(\tilde{Q})} \to 0.
\end{eqnarray*}
Finally, by Caccioppoli inequality
\[
\lim_{h\to+\infty}\sup_{y\in\tilde{K}}\|G_h(\cdot,y)-G(\cdot,y)\|_{H^1(Q)}=0.
\]
By the uniqueness of $G$ this holds for every subsequence of the original sequence $G_h$ so the claim follows.\qed

\begin{prop}\label{prop4.2} If for a positive $\ep_0$ and for some $k\in {1,\ldots, M}$ we have
\begin{equation}\label{7.1}
\left|S_k(y,z)\right|\leq r_0^{2-n}\ep_0\quad\mbox{for every }(y,z)\in K_0\times K_0,
\end{equation}
then
\[
    \left|S_k(\tilde{x},\tilde{x})\right|\leq
    C\left(\frac{\ep_0}{\ep_0+E}\right)^{\left(\tau^{(k+1)N_1}\delta_1^{k+1}\tau_r\right)^2}(\ep_0+E)r_0^{2-n}\left(\frac{r}{r_0}\right)^{2-n},
\]
where $\tilde{x}=P_{k+1}+r\nu_k$, $0<r<\frac{3r_0}{40}$, $\nu_k$ outer unit normal to $\der D_k$ at $P_{k+1}$ and
$\tau_r=\frac{\ln\left(\frac{3r_1-r}{3r_1-2r}\right)}{\ln\left(\frac{3r_1-r}{r_1}\right)}$.
\end{prop}

\textit{Proof of Proposition \ref{prop4.2}}.
Let us fix $z\in K_0$ and consider $v(y)=S_k(y,z)$. By Proposition \ref{prop4.1} we know that $v$ solves in $K$ equation
$\mbox{div}\left(\gamma^1\nabla v\right)=0$.
Moreover, by (\ref{3.20})
\[
\left|v(y)\right|\leq C E r_0^{1-\frac{n}{2}}\mbox{dist}(y,\Sigma_{k+1})^{1-\frac{n}{2}}\quad\mbox{for }y\in W_k,
\]
and by (\ref{7.1})
\[|v(y)|\leq r_0^{2-n}\ep_0\quad\mbox{for }y\in K_0.\]
Thus, we can apply  Proposition \ref{prop3.4} to $v$ getting for $0<r<\frac{3r_0}{40}$,
\begin{equation}\label{7.4}
    \left|S_k(\tilde{x},z)\right|\leq  \left(\frac{\ep_0}{E+\ep_0}\right)^\mu(E+\ep_0)\,
    r_0^{2-n}\left(\frac{r}{r_0}\right)^{(1-\frac{n}{2})(1-\tau_r)}.
\end{equation}
where we denoted by $\mu=\tau^{(k+1)N_1}\delta_1^{k+1}\tau_r$,
Now let us consider
\[
\tilde{v}(z)=S_k(\tilde{x},z)\quad\mbox{for } z\in K.
\]
which is solution of
\begin{equation}\label{7.5}
    \mbox{div}\left(\gamma^{(2)}\nabla \tilde{v}\right)=0\quad\mbox{in }K
\end{equation}
and satisfies
\begin{equation}\label{7.6}
    |\tilde{v}(z)|\leq CE\left(r \mbox{dist}(z,\Sigma_{k+1})\right)^{1-\frac{n}{2}}\quad\mbox{for }z\in K.
\end{equation}
By Proposition \ref{prop3.4} and taking into account (\ref{7.4}), (\ref{7.5}) and (\ref{7.6})
we have
\[
|\tilde{v}(\tilde{x})|\leq
C\left(\frac{\ep_0}{E+\ep_0}\right)^{\mu^2}(E+\ep_0)r^{2-n}.
\]
\qed

\textit{Proof of Theorem \ref{teo2.1}}
Let us  denote by $\ep:=\|\Lambda_{1}-\Lambda_{2}\|_{\mathcal{L}(H^{1/2},H^{-1/2})}$.
Let $\delta_k:=\| \gamma^{(1)}-\gamma^{(2)} \|_{L^\infty(W_k)}$ for $k=0,1,\ldots,M$.

Note that, for $y,z\in D_0$, due to the extension to the complex case of Alessandrini's identity (see formula (\ref{A3.4}) in the Appendix)
\[\int_{\om_0}(\gamma^{(1)}-\gamma^{(2)})\nabla G_1(\cdot,y)\nabla G_2(\cdot,z)=<(\Lambda_1-\Lambda_2)G_1(\cdot,y),G_2(\cdot,z)>\]
and, for $y,z\in K_0$,   Proposition \ref{prop3.3} yields
\begin{equation}\label{7.7}
    |S_{k-1}(y,z)|\leq C r_0^{2-n}(\ep+\delta_{k-1})
\end{equation}
Let $P_k\in\Sigma_k$ and $y_r=z_r=P_k+r\nu$ when $\nu$ is the outer normal vector to $\Sigma_k$ and $r\in (0,r_0/2)$.
Let us write
\begin{equation}\label{7.8}
  S_{k-1}(y_r,y_r)=I_1+I_2,
\end{equation}
where
\begin{eqnarray*}
  I_1 &=&  \int_{B_{\rho_0(P_k)\cap D_k}}(\gamma^{(1)}-\gamma^{(2)})\nabla G_1(\cdot,y_r)\nabla G_2(\cdot,y_r)\\
  I_2 &=&   \int_{U_{k-1}\setminus(B_{\rho_0(P_k)\cap D_k})}(\gamma^{(1)}-\gamma^{(2)})\nabla G_1(\cdot,y_r)\nabla G_2(\cdot,y_r).
\end{eqnarray*}
From Proposition \ref{prop3.3} we have that
\[
    |I_2|\leq C r_0^{2-n}E
\]
where $\rho_0=\frac{r_0}{6}$. On the other hand, again by Proposition \ref{prop3.3},
\begin{eqnarray*}\label{7.10}
    |I_1|&=&|\gamma^{(1)}-\gamma^{(2)}|\int_{B_{\rho_0(P_k)\cap D_k}}\!\!\!|\nabla G_1(\cdot,y_r)||\nabla G_2(\cdot,y_r)|\\
&\geq& C_1|\gamma^{(1)}_k-
\gamma^{(2)}_k|r^{2-n}
    -C_2 r_0^{2-n}E
\end{eqnarray*}

Now by (\ref{7.7}) and by Proposition \ref{prop4.2} we derive that, for $0<r<\frac{3r_0}{40}$,
\begin{equation}\label{7.11}
   |S_{k-1}(y_r,y_r)|\leq C
\left(\frac{\ep+\delta_{k-1}}{\ep+\delta_{k-1}+E}\right)^{(\tau^{kN_1}\delta_1^k\tau_r)^2}\!\!\!(\ep+\delta_{k-1}+E)r_0^{2-n}\left(\frac{r}{r_0}\right)^{2-n}
\end{equation}
Hence by (\ref{7.8})-(\ref{7.11}) we get
\[
 C_1|\gamma^{(1)}_k-
\gamma^{(2)}_k|r^{2-n}\leq C\left(r_0^{2-n}E+\left(\frac{\ep+\delta_{k-1}}{\ep+\delta_{k-1}+E}\right)^{(\tau^{kN_1}\delta_1^k\tau_r)^2}
 \!\!\!\!\!\!\!\!\!\!\!\!(\ep+\delta_{k-1}+E)r^{2-n}\right)
\]
and
\[
|\gamma^{(1)}_k-
\gamma^{(2)}_k|\leq C(\ep+\delta_{k-1}+E)\left[\left(\frac{\ep+\delta_{k-1}}{\ep+\delta_{k-1}+E}\right)^{(\tau^{kN_1}\delta_1^k\tau_r)^2}+\left(\frac{r}{r_0}\right)^{n-2}\right]
\]

By taking  $r=\left(\ln\left(\frac{\ep+\delta_{k-1}}{\ep+\delta_{k-1}+E}\right)\right)^{-1/4}$
we have
\begin{equation}\label{7.12}
    |\gamma^{(1)}_k-
\gamma^{(2)}_k|\leq C (\ep+\delta_{k-1}+E)\left(\ln\left(\frac{\ep+\delta_{k-1}}{\ep+\delta_{k-1}+E}\right)\right)^{-(n-2)/4}.
\end{equation}
Let
\[
\omega(t)=\left\{\begin{array}{ccc}
                   |\log t|^{-(n-2)/4} &\mbox{for}& 0<t\leq\frac{1}{e^n} \\
                   n^{-(n-2)/4} & \mbox{for}&t\geq\frac{1}{e^n}
                 \end{array}
\right.
\]
This function is increasing, concave, $\lim_{t\to 0}\omega(t)=0$ and the function $x\to x\omega(1/x)$ is increasing.
Inequality (\ref{7.12}) can be written as
\begin{equation}\label{dis}\delta_k\leq C(\ep+\delta_{k-1}+E)\omega\left(\frac{\ep+\delta_{k-1}}{\ep+\delta_{k-1}+E}\right).
\end{equation}
Notice that the above choice of $r$ is possible only if $\left(\ln\left(\frac{\ep+\delta_{k-1}}{\ep+\delta_{k-1}+E}\right)\right)^{-1/4}<\frac{3r_0}{40}$, but, if this is not the case, inequality (\ref{dis}) is obviously satisfied.

Since $\delta_0=0$ iterating (\ref{dis}) we obtain
\[
\delta_k+\ep\leq  (C+1)^k(E+\ep)\omega_k\left(\frac{\ep}{\ep+E}\right),
\]
where $\omega_k$ is the composition of $\omega$ with itself $k$ times.
Now we recall that $E=\delta_M$ and, hence,
\[
E+\ep\leq (C+1)^M(\ep+E)\omega_M\left(\frac{\ep}{\ep+E}\right)\]
where $C=C(\lambda,A,M,L)$.
Now, either $E\leq\ep$ and this proves Lipschitz stability, or $E>\ep$ and we can write
\[
E\leq 2(C+1)^M E\omega_M\left(\frac{\ep}{2E}\right),
\]
from which
\[\frac{1}{2(C+1)^M}\leq \omega_M\left(\frac{\ep}{2E}\right),\]
and
\[
\omega_M^{-1}\left(\frac{1}{2(C+1)^M}\right)\leq\frac{\ep}{2E},
\]
that is
\[
E\leq \frac{\ep}{2\omega_M^{-1}(\frac{1}{2(C+1)^M})}.
\]
Finally observing that $2\omega_M^{-1}(\frac{1}{2(C+1)^M})\geq 2\omega_N^{-1}(\frac{1}{2(C+1)^N})$ the claim follows.
\section{Concluding remarks}
\begin{rem}\label{remark2}
Observe that in our result we can replace the full Dirichlet to Neumann map with the local Dirichlet to Neumann map.
More precisely, let $\Sigma$ be an open portion of $\der\om$ containing a flat open subset.
Let $H^{1/2}_{co}(\Sigma)=\left\{\phi\in H^{1/2}(\der\om)\,:\,\mbox{supp }\phi\subset\Sigma\right\}$ and define the local
Dirichlet to Neumann map in the following way: for $\phi\in H^{1/2}_{co}(\Sigma)$ let
\[
<\Lambda_\gamma^\Sigma \phi,\psi>=\int_\om \gamma \nabla u\nabla \overline{v}\quad\mbox{for}\,\,\phi,\psi\in H^{1/2}_{co}(\Sigma)
\]
where $u$ solves equation $\mbox{div}(\gamma\nabla u)=0$ and $u=\phi$ on $\der\om$ and $v\in H^1(\om)$ such that $v=\psi$ on $\der\om$.

We observe that in our proof we apply the Dirichlet to Neumann map to functions whose support is contained in a
neighborhood of the flat portion of
$\der\om$. Hence in Theorem \ref{teo2.1} we can substitute (\ref{stabil})
by
\[
\|\gamma^{1}-\gamma^{2}\|_{L^\infty(\om)}\leq C\|\Lambda_1^\Sigma-\Lambda_2^\Sigma\|_{\mathfrak{L}(H^{1/2}_{co},H^{-1/2}_{co})},
\]
where
\[
\|\Lambda_\gamma^\Sigma\|_{\mathfrak{L}(H^{1/2}_{co},H^{-1/2}_{co})}=\sup\left\{<\Lambda_\gamma^\Sigma\phi,\psi>\,:\, \phi,\psi\in H^{1/2}_{co},\,\,
\|\phi\|_{H^{1/2}_{co}}=1,\|\psi\|_{H^{1/2}_{co}}=1\right\}
\]
\end{rem}
\begin{rem}\label{remark1}
We expect that  Lipschitz continuous dependence of the admittivities on the data still holds replacing the flatness condition on the
interfaces with $C^{1,\alpha}$ regularity.
In fact the key ingredients in our proof are the construction of a singular function $G(.,y)$ for $y\in K$, the unique
continuation estimates and the $C^{0,1}$ estimates in $K$ obtained by an application of the results contained in [LN].

In the case of $C^{1,\alpha}$ interfaces it is possible, proceeding similarly to what done in [AV], to make a $C^{1,\alpha}$ change of
variables which straightens locally the interface and to prove the existence of the Green's function in the new variables. Also unique
continuation estimates can be obtained in this case by means of three sphere inequality and finally $C^{\alpha}$ estimates in $K$ and local $C^{1,\alpha}$ estimates can be derived using the results in  [LN].
\end{rem}
\section{Appendix}
%
\subsection{Caccioppoli inequality}
%
For reader's convenience we state here Caccioppoli result for admittance equation. The proof is standard.
\begin{prop}\label{Caccio}
Let $u$ be a solution for
\begin{equation}\label{equazione}
    \mbox{div}\left(\gamma\nabla u\right)=0\quad\mbox{in}\quad \om,
\end{equation}
and let $B_{R}(x_0)\subset\om$. There exists a constant $C$ depending only on $\lambda$ such that
\begin{equation}\label{caccioppoli}
    \int_{B_\rho(x_0)}|\nabla u|^2\leq\frac{C}{(R-\rho)^2}\int_{B_R(x_0)}|u|^2,
\end{equation}
for every $\rho<R$.
\end{prop}

%
%
\subsection{Proof of Theorem \ref{teo3.3}.} Let us fix $x_0\in\Sigma_{j+1}$ such that $\mbox{dist}(x_0,P_{j+1})<\frac{r_0}{2}$ and
consider $B_R(x_0)$ with $R=\frac{r_0}{4}$.

By Proposition  (\ref{prop3.1}) (and the more general result in \cite{LN}), $D^{\beta^\prime}_{x^\prime} u$ is continuous
in $D_j$ and $D_{j+1}$ and belongs to  $C^1(D_j\cup \Sigma_{j+1})$ and $C^1(D_{j+1}\cup \Sigma_{j+1})$.

Observe that, since equation (\ref{equazione}) has constant coefficients
each $D^{\beta^\prime}_{x^\prime} u\in C(B_R(x_0))$ and it is a solution to the same equation for any multi-index $\beta$, hence
$u\in C^\infty(\overline{B}^+_R(x_0))$ and
$u\in C^\infty(\overline{B}^-_R(x_0))$.

We can apply Caccioppoli inequality (\ref{caccioppoli}) to
 $D^{\beta^\prime}_{x^\prime} u$ and get
\[
\int_{B_{\rho^{\prime\prime}}(x_0)}\left|\nabla(D^{\beta^\prime}_{x^\prime} u)\right|^2\leq
\frac{C}{(\rho^\prime-\rho^{\prime\prime})^2}\int_{B_{\rho^\prime}(x_0)}\left|D^{\beta^\prime}_{x^\prime} u\right|^2
\]
for $0<\rho^{\prime\prime}<\rho^\prime<R$ and $C=C(\lambda)$ of Proposition \ref{Caccio}.

Let $N_0$ be the maximum order of derivative that we want to estimate in $B_{\frac{R}{2}}(x_0)$ and define $\rho_k=R-\frac{kR}{2N_0}$, so that
$\rho_{k-1}-\rho_k=\frac{R}{2N_0}$, $k=0,\ldots,2N_0-1$. We have that
\[
  \int_{B_{\rho_1}(x_0)}|\nabla u|^2 \leq C\left(\frac{2N_0}{R}\right)^2\int_{B_{\rho_0}(x_0)}|u|^2,
\]
that is
\[
\sum_{k=1}^n\int_{B_{\rho_1}(x_0)}\left| \frac{\der u}{\der x_k}\right|^2\leq C\left(\frac{2N_0}{R}\right)^2\int_{B_{\rho_0}(x_0)}|u|^2.
\]
Since $\frac{\der u}{\der x_k}$ is solution to equation (\ref{equazione}) for every $k=1,\ldots,n-1$, Caccioppoli inequality holds and, consequently
\begin{eqnarray*}
  \int_{B_{\rho_2}(x_0)}\left|\nabla\left(\frac{\der u}{\der x_k} \right)\right|^2 &\leq& C\left(\frac{2N_0}{R}\right)^2\int_{B_{\rho_1}(x_0)}\left|\frac{\der u}{\der x_k} \right|^2 \\
   &\leq& \left(C\frac{2N_0}{R}\right)^2\int_{B_{\rho_0}(x_0)}|u|^2.
\end{eqnarray*}
Let now $\sigma_m=\sum_{|\beta|=m}\int_{B_{\rho_m}(x_0)}|D^{\beta^\prime}_{x^\prime} u|^2$ and assume that
\[
\sigma_m\leq \left(C\left(\frac{2N_0}{R}\right)^2\right)^m\int_{B_{\rho_0}(x_0)}|u|^2.
\]
We proceed by induction and consider
\begin{eqnarray*}
&&\sigma_{m+1}=\sum_{|\beta|=m+1}\int_{B_{\rho_{m+1}(x_0)}}|D^{\beta^\prime}_{x^\prime} u|^2=
\sum_{|\beta|=m}\sum_{k=1}^{n-1}\int_{B_{\rho_{m+1}}(x_0)}\left|\frac{\der}{\der x_k}D^{\beta^\prime}_{x^\prime} u\right|^2
\\&&\leq \frac{C}{(\rho_m-\rho_{m-1})^2}\sum_{|\beta|=m}\int_{B_{\rho_m}(x_0)}\!\!\!|D^{\beta^\prime}_{x^\prime} u|^2=
C\left(\frac{2N_0}{R}\right)^2\sum_{|\beta|=m}\int_{B_{\rho_m}(x_0)}\!\!\!|D^{\beta^\prime}_{x^\prime} u|^2\\
&&\leq\left(C\left(\frac{2N_0}{R}\right)^2\right)^{m+1}\int_{B_{\rho_0}(x_0)}|u|^2.
\end{eqnarray*}
Hence, we have proved in particular that
\begin{equation}\label{A3.18}
    \sum_{|\beta|=N_0}\int_{B_{\frac{R}{2}}(x_0)}|D^{\beta^\prime}_{x^\prime} u|^2\leq \left(C\left(\frac{2N_0}{R}\right)^2\right)^{N_0}
    \int_{B_R(x_0)}|u|^2.
\end{equation}
Finally, applying  Proposition \ref{prop3.1}, we have that
\[
\|D^{\beta^\prime}_{x^\prime} u\|_{L^\infty(B_{\frac{R}{4}}(x_0))}\leq
CR^{-\frac{n}{2}}\|D^{\beta^\prime}_{x^\prime}u\|_{L^2(B_{\frac{R}{2}(x_0)})}.
\]
Hence, by (\ref{A3.18}) and recalling that $N_0=|\beta|$ we get
\[
\|D^{\beta^\prime}_{x^\prime} u\|_{L^\infty(B_{\frac{R}{4}}(x_0))}\leq C\left(\frac{1}{R^n}\int_{B_{\frac{R}{2}}(x_0)}|u|^2\right)^{1/2}
\left(C\left(\frac{2|\beta|}{R}\right)^2\right)^{\frac{|\beta|}{2}}.
\]
Observing that $|\beta|^\beta\leq n^{|\beta|}\beta!e^{|\beta|}$ one gets
\begin{equation}\label{A3.19}
    \|D^{\beta^\prime}_{x^\prime} u\|_{L^\infty(B_{\frac{R}{4}}(x_0))}\leq M \beta! \left(\frac{C_1}{R}\right)^{|\beta|}
\end{equation}
where $M=C\left(\frac{1}{R^n}\int_{B_{\frac{R}{2}}(x_0)}|u|^2\right)^{\frac{1}{2}}$ and $C_1=e n8\lambda^2$.

In particular, from (\ref{A3.19}) we derive
\begin{equation}\label{A3.20}
    \|D^{\beta^\prime}_{x^\prime} u(\cdot,0)\|_{L^\infty(B^{\prime}_{\frac{R}{4}}(x_0))}\leq M\beta!\left(\frac{C_1}{R}\right)^{|\beta|}
\end{equation}
which implies analyticity of $\Re u(x^\prime,0)$ and $\Im u(x^\prime,0)$.

Let $\phi(x^\prime)=\Re u(x^\prime,0)$ and set $U(x)=:\Re u(x)-\phi(x^\prime)$. Then $U$ satisfies
\[
    \left\{
\begin{array}{rcl}
  \Delta U & = &-\Delta\phi :=F\quad\mbox{in}\quad B_R^-(x_0) \\
 U(x^\prime,0) &= & 0\quad\mbox{on}\quad B_R^{\prime}(x_0).
\end{array}
    \right.
\]
Observe now that, by standard regularity estimates (cfr. for example \cite[Corollary 2.36]{GT}) we have that
\[
\left\|\nabla U\right\|_{L^\infty(B_{\frac{R}{8}}^-(x_0))}\leq
C\left(\|U\|_{L^\infty(B_{\frac{R}{4}}^-(x_0))}+
\|F\|_{L^\infty(B_{\frac{R}{4}}^-(x_0))}\right).
\]
Analogously, the function $W:=D^{\beta^\prime}_{x^\prime} U$, solves the problem
\[
\left\{
\begin{array}{rcl}
\Delta W&=& -\Delta (D^{\beta^\prime}_{x^\prime}\phi)\quad\mbox{in}\quad B_R^-(x_0)\\
W&=&0\quad\mbox{on}\quad B_R^\prime(x_0),\end{array}
\right.
\]
so that also
\[
\|\nabla W\|_{L^\infty(B_{\frac{R}{8}}^-(x_0))}\leq C\left(\|W\|_{L^\infty(B_{\frac{R}{4}}^-(x_0))}
+\|\Delta D^{\beta^\prime}_{x^\prime}\phi\|_{L^\infty(B_{\frac{R}{4}}^-(x_0))}\right)
\]
Applying estimates (\ref{A3.19}) and (\ref{A3.20}) we then get that
\[
\left\|D^{\beta^\prime}_{x^\prime} \frac{\der }{\der x_n} U\right\|_{L^\infty(B^\prime_{\frac{R}{8}}(x_0))}\leq\frac{C}{R^4}M\beta!
    \left(\frac{C_1}{R}\right)^{|\beta|}.
\]
Hence, also,
\begin{equation}\label{A3.22}
\|D^{\beta^\prime}_{x^\prime} \frac{\der }{\der x_n}\Re u\|_{L^\infty(B^\prime_{\frac{R}{8}}(x_0))}\leq\frac{C}{R^4}M\beta!
    \left(\frac{C_1}{R}\right)^{|\beta|},
\end{equation}
which implies analyticity of $ \frac{\der }{\der x_n}\Re u(x',0)$.
A similar estimate can be proved for $\Im u$. By the fact that  $\Re u$ and $\Im u$ are harmonic in $B^+$ and in $B^-$ separately,
and  from (\ref{A3.20}) and (\ref{A3.22}) we have that,
\[
\left|\Re u(x)\right|=\left|\sum_\beta \frac{D^{\beta}_{x}\Re u(x_0)}{\beta!}(x-x_0)^\beta\right|\leq
\frac{C}{R^4}M\sum_{\beta}\left(\frac{C_1}{R}\right)^{|\beta|}|x-x_0|^{|\beta|}
\]
which is convergent for
\[|x-x_0|=\tilde{R}=\frac{R}{2C_1}.\]
Hence $\Re u$ can be extended analytically in a neighborhood of $x_0$ and the same can be proved for $\Im u$.
Repeating the same argument for all points  $x\in\Sigma_{j+1}$ such that $\mbox{dist}(x,P_{j+1})<\frac{r_0}{2}$ and
choosing $R=\frac{r_0}{4}$ we have proved that $\Re u$ can be extended analytically to the set
$E_{j+1}^{(C_1)}$.

The same holds true for $\Im u$, hence $u_{D_j}$ can be extended analytically to $D_j\cup E_{j+1}^{(C_1)}$
and the extension $\tilde{u}$ satisfies (\ref{3.32}).


%
%
\subsection{A generalization of Alessandrini's identity}
Alessandrini's identity holds for solutions to conductivity equation with real valued coefficients. For sake of completeness we show here that is can be generalized to the case of complex valued coefficient.

Let $u_1$ and $u_2$ be the solutions to
\begin{equation}\label{A3.3}
    \mbox{div}\left(\gamma^{(k)}\nabla u_k\right)=0\quad\mbox{in}\quad \om,
\end{equation}
for $k=1,2$ respectively and let us consider the Dirichlet-to-Neumann maps $\Lambda_{\gamma^{(k)}}$ that, from now on,
we will denote by $\Lambda_k$.

From (\ref{A3.3}) we have
\[
\int_\om\gamma^{(1)}\nabla u_1\cdot\nabla\ov{v}=<\Lambda_1 u_1, v>,\quad\forall v\in H^1(\om).
\]
If we take $v=\ov{u}_2$, we derive
\[
\int_\om\gamma^{(1)}\nabla u_1\cdot\nabla u_2=<\Lambda_1 u_1, \ov{u}_2>,
\]
and, analogously,
\[
\int_\om\gamma^{(2)}\nabla u_2\cdot\nabla u_1=<\Lambda_2 u_2, \ov{u}_1>.
\]

Let us show that $<\Lambda_1 u_1, \ov{u}_2>=<\Lambda_1 u_2, \ov{u}_1>$. Let $w$ be solution to
\[
    \left\{\begin{array}{rcl}
             \mbox{div}\left(\gamma^{(1)}\nabla w\right)&= & 0\quad\mbox{in}\quad\om \\
             w&= &u_2  \quad\mbox{on}\quad\der\om,\\
           \end{array}
    \right.
\]
then,
\[
<\Lambda_1 u_2, \ov{u}_1>=\int_\om\gamma^{(1)}\nabla w\cdot\nabla u_1= <\Lambda_1 u_1, \ov{u}_2>.
\]
Hence
\begin{equation}\label{A3.4}
    \int_\om (\gamma^{(1)}-\gamma^{(2)})\nabla u_1\cdot\nabla u_2=<\left(\Lambda_1 -\Lambda_2\right)u_2, \ov{u}_1>.
\end{equation}
\section*{Acknowledgements}
\noindent
We would like to thank Micol Amar, Daniele Andreucci, Paolo Bisegna and Roberto Gianni for having pointed out the importance of this problem in connection to the study of a model of conduction in biological tissues.

\end{document}